\documentclass{article}[12pt]
\usepackage{amssymb}
\usepackage{eurosym}
\usepackage{epsfig}
\usepackage{amsmath}
\usepackage{amsfonts}
\usepackage{pgf,tikz}
\usepackage{enumerate}
\usepackage{cite}
\usepackage{xcolor, colortbl}

\usepackage[colorlinks=true, linkcolor=black, citecolor=black, urlcolor=black]{hyperref}

\newcommand{\R}{{{\Bbb R}}}

\newtheorem{theorem}{\sc Theorem}[section]
\newtheorem{proposition}{\sc Proposition}[section]
\newtheorem{lemma}{\sc Lemma}[section]

\newtheorem{remark}{\sc Remark}[section]
\newtheorem{corollary}{\sc Corollary}[section]
\newtheorem{example}{\sc Example}[section]

\def\qed{\hbox to 0pt{}\hfill$\rlap{$\sqcap$}\sqcup$\medbreak}

\title{On Krasnosel'ski\u{\i}--Precup fixed point theorem and eigenvalue criteria for existence of positive solutions}

\author{Laura Mª Fern\'andez--Pardo and Jorge Rodr\'iguez--L\'opez} 

\date{}
\begin{document}
 \maketitle

\begin{center}  {\small CITMAga \& Departamento de Estat\'{\i}stica, An\'alise Matem\'atica e Optimizaci\'on, \\ Universidade de Santiago de Compostela, \\ 15782, Facultade de Matem\'aticas, Campus Vida, Santiago, Spain.\\  Email: lauramaria.fernandez.pardo@rai.usc.es; jorgerodriguez.lopez@usc.es}
\end{center}

\medbreak

\noindent {\it Abstract.} We provide an alternative approach to Krasnosel'ski\u{\i}--Precup fixed point theorem for operator systems based on the Leray-Schauder fixed point index in cones. Our focus is on the case of operators whose components are entirely of compressive type. The abstract technique is applied to a system of second--order differential equations providing a coexistence positive solution by means of an eigenvalue type criterion.  

\medbreak


\medbreak


\section{Introduction}

In this paper, we are concerned with the so called \textit{vector version of Krasnosel'ski\u{\i} fixed point theorem} or \textit{Krasnosel'ski\u{\i}--Precup fixed point theorem}, that is, an existence result for systems of operator equations of the form
\begin{equation*}
	\left\{\begin{array}{l} x_1=T_1(x_1,x_2), \\ x_2=T_2(x_1,x_2), \end{array} \right.
\end{equation*}
which imposes compression--expansion conditions on each component of the system. One of the main features of this result due to Precup (see \cite{PrecupFPT,PrecupSDC}) is that all the components of the obtained fixed point are non-trivial, i.e., it is a \textit{coexistence} fixed point.  In fact, each component of the fixed point $(x_1,x_2)$ is localized independently in norm between two positive constants in the following way: $r_1\leq\left\|x_1\right\|\leq R_1$ and $r_2\leq\left\|x_2\right\|\leq R_2$.

Our approach relies on the fixed point index for compact maps in cones and it is intended to complement the original version due to Precup. It continues the strategy initiated in \cite{JRL} but combined with the idea introduced in \cite{LFP_JRL}, from which we derive a refinement in the compressive case. The key ingredient is the choice of an appropriate retraction when extending the operator $T$, ensuring the absence of new fixed points under this less restrictive situation and thus yielding a well-defined index for the computations, as explained below.  Note that the calculation of the fixed point index of the corresponding operator $T=(T_1,T_2)$ provides useful information for deriving multiplicity criteria. 

As an application we study the existence of positive solutions to a system of second--order ODEs subject to Robin type boundary conditions. Our criterion, which is inspired by the previous work of Webb, Infante and Lan \cite{WI,WL,W}, involves the asymptotic behavior of the nonlinearities at $0$ and infinity in relation to the principal eigenvalue of the corresponding linear Hammerstein type operator associated to the differential problem. To our best knowledge, it is the first time that the vector approach is combined with this sharp technique based on eigenvalue type conditions, yielding solutions which are non-trivial in each of its components. Furthermore, the reasoning for establishing this existence result cannot be applied under the assumptions in \cite{JRL}, but is valid under the weaker conditions presented here.

The paper is structured as follows. In  Section 2 we review the vectorial versions of the Krasnosel’skiĭ fixed point theorem, we recall some basic properties of the Leray-Schauder fixed point index in cones and compute it in the compression case under both order relation conditions and homotopy conditions. Section 3 addresses the application to the existence of solutions with positive components for the second-order problem previously described.

\section{On Krasnosel'ski\u{\i}--Precup fixed point theorem in cones}

\subsection{Preliminaries on Krasnosel'ski\u{\i} type fixed point theorems in Cartesian products}

A closed convex subset $K$ of a normed linear space $X$ is said to be a \textit{cone} if $\lambda\,u\in K$ for every $u\in K$ and for all $\lambda\geq 0$, and $K\cap (-K)=\{0\}$. A cone $K$ induces a partial order in $X$ given by $x\preceq y$ if and only if $y-x\in K$. In addition, we will use the notation $x\prec y$ in case that $y-x\in K\setminus\{0\}$.

Let us recall now the fixed point theorem for operator systems (or operators defined in the Cartesian product of cones) proposed by Precup \cite{PrecupSDC}, which we denominate \textit{Krasnosel'ski\u{\i}--Precup fixed point theorem in cones} since the compression--expansion conditions in the line of Krasnosel'ski\u{\i} are imposed in a component-wise manner. In the sequel, consider two cones $K_1$ and $K_2$ of the normed linear spaces $X$ and $Y$, and so $K:=K_{1}\times K_{2}$ is a cone of $X\times Y$. For $r,R\in \mathbb{R}^{2}$, $r=(r_{1},r_{2})$, $R=(R_{1},R_{2})$, with $0<r_i<R_i$ $(i=1,2)$, we denote
\[ K_{r,R}:=\{x=\left( x_{1},x_{2}\right) \in K\,:\,r_{i}< \left\Vert x_{i}\right\Vert < R_{i}\quad \text{for }i=1,2\}=(K_1)_{r_1,R_1}\times (K_2)_{r_2,R_2} \]
and
\[ \overline{K}_{r,R}:=\{x=\left( x_{1},x_{2}\right) \in K\,:\,r_{i}\leq \left\Vert x_{i}\right\Vert \leq R_{i}\quad \text{for }i=1,2\}=(\overline{K}_1)_{r_1,R_1}\times (\overline{K}_2)_{r_2,R_2}.\]
The result due to Precup is the following.

\begin{theorem}[Krasnosel'ski\u{\i}-Precup]\label{th_KP}
	Let $r,R\in\R^2$, $r=(r_{1},r_{2})$, $R=(R_{1},R_{2})$, with $0<r_i<R_i$ $(i=1,2)$.	
	Assume that $T=(T_1,T_2):\overline{K}_{r,R}\rightarrow K$ is a compact map and for each $i\in\{1,2\}$ there exists $h_i\in K_i\setminus\{0\}$ such that one of the following conditions is satisfied in $\overline{K}_{r,R}$:
	\begin{enumerate}
		\item[$(a)$] $T_i x+\mu\,h_i\neq x_i$ if $\left\|x_i\right\|=r_i$ and $\mu>0$, and $T_i x\neq \lambda\, x_i$ if $\left\|x_i\right\|=R_i$ and $\lambda>1$;
		\item[$(b)$] $T_i x\neq \lambda\, x_i$ if $\left\|x_i\right\|=r_i$ and $\lambda>1$, and $T_i x+\mu\,h_i\neq x_i$ if $\left\|x_i\right\|=R_i$ and $\mu>0$. 
	\end{enumerate}
	Then $T$ has at least a fixed point $\bar{x}=(\bar{x}_1,\bar{x}_2)\in K$ with $r_i\leq\left\|\bar{x}_i\right\|\leq R_i$ $(i=1,2)$.
\end{theorem}

\begin{remark}
	Conditions $(a)$ and $(b)$ can be respectively replaced by the following stronger assumptions:
	\begin{enumerate}
	\item[$(a^{\star})$] $T_i x\nprec x_i$ if $\left\|x_i\right\|=r_i$  and $T_i x\nsucc x_i$ if $\left\|x_i\right\|=R_i$;
	\item[$(b^{\star})$] $T_i x\nsucc x_i$ if $\left\|x_i\right\|=r_i$ and $T_i x\nprec x_i$ if $\left\|x_i\right\|=R_i$. 
	\end{enumerate}
	Theorem \ref{th_KP} with conditions $(a)$ and $(b)$ replaced by $(a^{\star})$ and $(b^{\star})$ is the main result of the paper by Precup \cite{PrecupFPT}. It can be seen as a particular case of those later obtained by the same author in \cite{PrecupSDC}.  
\end{remark}

The proof of Theorem \ref{th_KP} presented in \cite{PrecupFPT} is based on Schauder fixed point theorem in case that condition $(a)$ holds for both $i=1$ and $i=2$. Then the remaining cases are reduced to the solved one by means of a trick which transforms operators satisfying $(b)$ into operators fulfilling $(a)$. A different approach, based on the fixed point index in cones, has been proposed in \cite{JRL}. As already explained in \cite{JRL}, the computation of the index has interest itself and provides more information than merely the existence of a fixed point (for instance, in order to extend the result to more general contexts or to derive multiplicity criteria). Note that the result presented in \cite{JRL} requires slightly stronger conditions on the operator when computing the fixed point index. It can be summarized as follows.

\begin{theorem}\label{th_KP_strong}
	Let $r,R\in\R^2$, $r=(r_{1},r_{2})$, $R=(R_{1},R_{2})$, with $0<r_i<R_i$ $(i=1,2)$.	
	Assume that $T=(T_1,T_2):\overline{K}_{r,R}\rightarrow K$ is a compact map and for each $i\in\{1,2\}$ there exists $h_i\in K_i\setminus\{0\}$ such that one of the following conditions is satisfied in $\overline{K}_{r,R}$:
	\begin{enumerate}
		\item[$(a^{\dagger})$] $T_i x+\mu\,h_i\neq x_i$ if $\left\|x_i\right\|=r_i$ and $\mu\geq 0$, and $T_i x\neq \lambda\, x_i$ if $\left\|x_i\right\|=R_i$ and $\lambda\geq 1$;
		\item[$(b^{\dagger})$] $T_i x\neq \lambda\, x_i$ if $\left\|x_i\right\|=r_i$ and $\lambda\geq 1$, and $T_i x+\mu\,h_i\neq x_i$ if $\left\|x_i\right\|=R_i$ and $\mu\geq 0$. 
	\end{enumerate}
	Then \[i_{K}(T,K_{r,R})=(-1)^k \]
	where $k\in\{0,1,2\}$ denotes the number of indexes for which condition $(b^{\dagger})$ holds.
	
	In particular, $T$ has at least a fixed point $\bar{x}=(\bar{x}_1,\bar{x}_2)\in K$ with $r_i<\left\|\bar{x}_i\right\|< R_i$ $(i=1,2)$.
\end{theorem}  

\begin{remark}
	Observe that conditions $(a^{\dagger})$ - $(b^{\dagger})$ are more stringent than conditions $(a)$ - $(b)$ together with the fact that $T$ is fixed point free on the boundary of $K_{r,R}$. Indeed, $(a^{\dagger})$ - $(b^{\dagger})$ imply that
	\[T_i x\neq x_i \ \text{ if } \left\|x_i\right\|=r_i \text{ or if } \left\|x_i\right\|=R_i, \quad i=1,2, \]
	which, in turn, ensures that $T=(T_1,T_2)$ has no fixed points on $\partial_{K} K_{r,R}$. The reverse implication is clearly false.
\end{remark}

In view of the previous remark, we wonder if it is possible to compute the \textit{fixed point index} of the operator $T$ over the set $K_{r,R}$ under the assumptions $(a)$ - $(b)$ (or $(a^{\star})$ - $(b^{\star})$) instead of $(a^{\dagger})$ - $(b^{\dagger})$. This is the main goal of the remaining part of the present section. We success in some particular (but meaningful) cases, as shown below.

\subsection{Preliminaries on the fixed point index in cones}

For the convenience of the reader, let us recall some concepts and properties related to the fixed point index in cones. 
If $U$ is a relatively open bounded subset of a cone $K$ of a normed space $X$ and $N:\overline{U}\rightarrow K$ is a compact map without fixed points on the boundary of $U$ (denoted by $\partial_K\,U$), the \textit{fixed point index} of $N$ on $U$ with respect to the cone $K$, $i_K(N,U)$, is well-defined. For additional details, we refer to \cite{amann,Granas,guolak}. 


\begin{proposition}\label{prop_index}
	Let $K$ be a cone of a normed space, $U\subset K$ be a bounded relatively open set and $N:\overline{U}\rightarrow K$ be a compact map such that $N$ has no fixed points on the boundary of $U$. Then the fixed point index of $N$ on the set $U$ with respect to $K$, $i_{K}(N,U)$, has the following properties:
	\begin{enumerate}
		\item (Additivity) Let $U$ be the disjoint union of two open sets $U_1$ and $U_2$. If $0\not\in(I-N)(\overline{U}\setminus(U_1\cup U_2))$, then \[i_{K}(N,U)=i_{K}(N,U_1)+i_{K}(N,U_2).\]
		\item (Existence) If $i_{K}(N,U)\neq 0$, then there exists $x\in U$ such that $Nx=x$.
		\item (Homotopy invariance) If $H:\overline{U}\times[0,1]\rightarrow K$ is a compact homotopy and $0\not\in(I-H)(\partial_K\,U\times[0,1])$, then
		\[i_{K}(H(\cdot,0),U)=i_{K}(H(\cdot,1),U).\]
		\item (Normalization) If $N$ is a constant map with $Nx=\bar{x}$ for every $x\in\overline{U}$, then
		\[i_{K}(N,U)=\left\{\begin{array}{ll} 1, & \text{ if } \bar{x}\in U, \\ 0, & \text{ if } \bar{x}\not\in\overline{U}. \end{array} \right. \]
	\end{enumerate}
\end{proposition}

The following computations of the fixed point index will be also useful in our reasoning. Their proofs can be found, for instance, in \cite[Lemma 2.3.1 and 2.3.2]{guolak}.

\begin{proposition}\label{prop_ind01}
	Let $K$ be a cone, $U\subset K$ be a bounded relatively open set such that $0\in U$ and $N:\overline{U}\rightarrow K$ be a compact map without fixed points on $\partial_{K}\,U$.
	\begin{enumerate}[$(i)$]
		\item If $Nx\neq \lambda\, x$ for all $x\in \partial_{K}\, U$ and all $\lambda> 1$, then $i_{K}(N,U)=1$.
		\item If there exists $h\in K\setminus\{0\}$ such that $Nx+\mu\, h\neq x$ for all $x\in \partial_{K}\, U$ and all $\mu> 0$, then $i_{K}(N,U)=0$.
	\end{enumerate}	
\end{proposition}

In case of operators defined in the Cartesian product of normed spaces, the following computation of the fixed point index is useful. Its proof can be found in \cite{PreRod,JRL}.

\begin{proposition}\label{prop_ind_sys}
	Let $U\times V$ be a bounded relatively open subset of the cone $K= K_1\times K_2$ in the normed spaces product $X=X_1\times X_2$, such that $0\in U$. Assume that $N=(N_1,N_2):\overline{U\times V}\rightarrow K$ is a compact mapping and there exists $h\in K_2\setminus\{0\}$ satisfying the following conditions:
	\begin{enumerate}[$(i)$]
		\item $N_1 x \neq \lambda\, x_1$ for all $x_1\in \partial_{K_1}U$, $x_2\in\overline{V}$ and all $\lambda>1$;
		\item $N_2 x+ \mu\, h\neq x_2$ for all $x_1\in\overline{U}$, $x_2\in\partial_{K_2} V$ and all $\mu>0$.
	\end{enumerate}
	If $N$ has no fixed points on $\partial_K\, (U\times V)$, then $i_{K}(N,U\times V)=0$.
\end{proposition}

\subsection{Krasnosel'ski\u{\i}--Precup fixed point theorem with order type conditions}

In this subsection, we focus on the computation of the fixed point index under compression type assumptions given by the order relation provided by the cone. More specifically, we will prove the following result.

\begin{theorem}\label{th_Korder}
	Let $r,R\in\R^2$, $r=(r_{1},r_{2})$, $R=(R_{1},R_{2})$, with $0<r_i<R_i$ $(i=1,2)$.	
	Assume that $T=(T_1,T_2):\overline{K}_{r,R}\rightarrow K$ is a compact map and for each $i\in\{1,2\}$ the following condition is satisfied in $\overline{K}_{r,R}$:
	\begin{enumerate}
			\item[$(a^{\star})$] $T_i x\nprec x_i$ if $\left\|x_i\right\|=r_i$  and $T_i x\nsucc x_i$ if $\left\|x_i\right\|=R_i$.
	\end{enumerate}
	If $T$ has no fixed points on the boundary of $K_{r,R}$, then $i_{K}(T,K_{r,R})=1$. 
\end{theorem}

To do so, adapting the reasoning in \cite{JRL}, we need to extend the definition of $T$ to a larger set, namely, $\overline{K}_{R}:=\{x=\left( x_{1},x_{2}\right) \in K\,:\, \left\Vert x_{i}\right\Vert \leq R_{i}\quad \text{for }i=1,2\}$. Such extension is based on the fact that $\overline{K}_{r,R}$ is a retract of $\overline{K}_R$. A retraction can be defined as follows:
\begin{equation}\label{eq_r1}
	\rho(x)=\left(\rho_1(x_1),\rho_2(x_2)\right), \quad \rho_i(x_i)=\left\{\begin{array}{ll} r_i\displaystyle\frac{x_i+(r_i-\left\|x_i\right\|) h_i}{\left\|x_i+(r_i-\left\|x_i\right\|) h_i \right\|}, & \quad \text{if } \left\|x_i\right\|<r_i, \\[0.2cm] x_i,  & \quad \text{if } r_i\leq\left\|x_i\right\|\leq R_i,  \end{array} \right.
\end{equation}
where $h_i\in K_i$, with $\left\|h_i\right\|=1$, is fixed. Using this retraction $\rho$, a map $T$ under the assumptions of Theorem~\ref{th_Korder} can be extended to the set $\overline{K}_R$ in the following way:
\begin{equation}\label{eq_N}
	N=(N_1,N_2):\overline{K}_R\rightarrow K, \quad N:=T\circ\rho.
\end{equation}
The proof of Theorem \ref{th_Korder} is finally based on the computation of the fixed point index of the map $N$ on different open sets. In order to ensure that the index is well-defined, it is necessary that the corresponding operator $N$ have no fixed points on the boundaries of these open sets. This is the case when $T$ satisfies condition $(a^{\star})$ and has not fixed points on $\partial_K\,K_{r,R}$, as shown by the following result.

\begin{lemma}\label{lem_noFP}
	Let $T=(T_1,T_2):\overline{K}_{r,R}\rightarrow K$ be a compact map such that for each $i\in\{1,2\}$ the following condition is satisfied in $\overline{K}_{r,R}$:
	\[T_i x\nprec x_i \ \text{ if } \ \left\|x_i\right\|=r_i. \]
	
	If $T$ has no fixed point on the boundary of $K_{r,R}$, then its extension $N$ defined as \eqref{eq_N} has also no fixed point $x=(x_1,x_2)\in\overline{K}_{R}$ such that $\left\|x_i\right\|=r_i$ or $\left\|x_i\right\|=R_i$ for $i\in\{1,2\}$.	
\end{lemma} 

\noindent
{\bf Proof.} Let us suppose that $N$ has a fixed point $(x_1,x_2)\in \overline{K}_{R}$ such that $\left\|x_1\right\|\in \{r_1,R_1\}$. Since $N=T$ in $\overline{K}_{r,R}$ and $T$ has no fixed points on the boundary of $\overline{K}_{r,R}$, it follows that $\left\|x_2\right\|<r_2$. Then
\[(x_1,x_2)=N(x_1,x_2)=\left(T_1(x_1,\rho_2(x_2)),T_2(x_1,\rho_2(x_2))\right).\]
In particular, $x_2=T_2(x_1,\rho_2(x_2))$ and, moreover, by the definition of $\rho_2$, we have $x_2\prec \rho_2(x_2)$. Indeed, $\rho_2(x_2)-x_2\in K_2\setminus\{0\}$ since
\[\rho_2(x_2)-x_2=\displaystyle\left(\frac{r_2}{\left\|x_2+(r_2-\left\|x_2\right\|) h_2\right\|}-1\right)x_2+\frac{r_2-\left\|x_2\right\|}{\left\|x_2+(r_2-\left\|x_2\right\|) h_2\right\|}h_2, \]
where the first addend belongs to $K_2$ because $\left\|x_2+(r_2-\left\|x_2\right\|) h_2\right\|\leq\left\|x_2\right\|+(r_2-\left\|x_2\right\|) \left\|h_2\right\|=r_2$ and the second one belongs to $K_2\setminus\{0\}$ as a consequence of $h_2\in K_2\setminus\{0\}$.

In conclusion, $T_2(x_1,\rho_2(x_2)) \prec \rho_2(x_2)$ and $\left\|\rho_2(x_2)\right\|=r_2$, a contradiction with the assumption.
\qed

Now we are in a position to demonstrate Theorem \ref{th_Korder}, whose proof relies on the properties of the fixed point index. The proposed computation of the index follows previous ideas employed in \cite{JRL}, but we include them here again for the convenience of the reader.

\medskip

\noindent
{\bf Proof of Theorem \ref{th_Korder}.} First, note that $(a^{\star})$ implies that for each $i\in\{1,2\}$ the following condition is satisfied in $\overline{K}_{R}$:
\[N_i x\nprec x_i \ \text{ if } \ \left\|x_i\right\|=r_i \ \text{ and } \ N_i x\nsucc x_i \ \text{ if } \ \left\|x_i\right\|=R_i. \]
As a consequence (see for instance \cite{guolak}), we have that for every $h_i\in K_i\setminus\{0\}$,
\[N_i x+\mu\,h_i\neq x_i \ \text{ if } \ \left\|x_i\right\|=r_i \ \text{ and } \ \mu> 0, \text{ and } \ N_i x\neq \lambda\, x_i \ \text{ if } \ \left\|x_i\right\|=R_i \ \text{ and } \ \lambda> 1. \]
Moreover, notice that Lemma \ref{lem_noFP} guarantees that the fixed point index of the map $N$ is well-defined on the open sets $K_r$, $K_R$, $(K_1)_{r_1}\times (K_2)_{R_2}$ and $(K_1)_{R_1}\times (K_2)_{r_2}$.
 
Hence, by Proposition \ref{prop_ind01} we have
\[i_{K}(N,K_{R})=1 \quad \text{and} \quad i_{K}(N,K_{r})=0. \]
In addition, Proposition \ref{prop_ind_sys} ensures that
\[i_{K}(N,(K_1)_{r_1}\times (K_2)_{R_2})=0=i_{K}(N,(K_1)_{R_1}\times(K_2)_{r_2}). \]
Applying the additivity property of the index twice,
\[i_{K}(N,(K_1)_{r_1,R_1}\times(K_2)_{r_2})=i_{K}(N,(K_1)_{R_1}\times(K_2)_{r_2})-i_{K}(N,K_r)=0 \]
and thus 
\[i_K(N,K_{r,R})=i_{K}(N,K_{R})-i_{K}(N,(K_1)_{r_1,R_1}\times(K_2)_{r_2})-i_{K}(N,(K_1)_{r_1}\times(K_2)_{R_2})=1. \]
In conclusion, $i_K(T,K_{r,R})=1$, since $N=T$ on the set $\overline{K}_{r,R}$.
\qed 

\begin{corollary}
	Let $r,R\in\R^2$, $r=(r_{1},r_{2})$, $R=(R_{1},R_{2})$, with $0<r_i<R_i$ $(i=1,2)$.	
	Assume that $T=(T_1,T_2):\overline{K}_{r,R}\rightarrow K$ is a compact map and for each $i\in\{1,2\}$  condition $(a^{\star})$ is satisfied in $\overline{K}_{r,R}$.
	
	Then $T$ has at least one fixed point in $\overline{K}_{r,R}$.
\end{corollary}

\noindent
{\bf Proof.} Either $T$ has a fixed point on the boundary of $\overline{K}_{r,R}$ (and the proof is finished) or, if not, Theorem \ref{th_Korder} is applicable and ensures that $i_K(T,K_{r,R})=1\neq 0$. Then the conclusion is deduced from the existence property of the index.
\qed  

\begin{remark}
	Theorem \ref{th_Korder} remains valid (under the same proof!) for operators defined in the Cartesian product of the set difference of two nested bounded, neighborhoods of the origin of the form \[(\overline{\Omega}_1\setminus (K_1)_{r_1})\times (\overline{\Omega}_2\setminus (K_2)_{r_2}),\] being $\Omega_1$ and $\Omega_2$ relatively open sets.
	In that case, condition $(a^{\star})$ should be replaced by the following:
	\begin{enumerate}
		\item[$(A^{\star})$]  $T_i x\nprec x_i$ if $\left\|x_i\right\|=r_i$  and $T_i x\nsucc x_i$ if $x_i\in\partial_{K_i}\Omega_i$.
	\end{enumerate}
	
\end{remark}

\subsection{Krasnosel'ski\u{\i}--Precup fixed point theorem with homotopy type conditions}

Now, we are concerned with cone--compression conditions given by means of homotopy type assumptions. Our aim is to provide a computation of the fixed point index in the line of Theorem~\ref{th_Korder}, but with assumption $(a)$ instead of $(a^{\star})$. Then, an analogous of Lemma \ref{lem_noFP} needs be proven and the definition of an \textit{adequate} retraction becomes crucial at this step. 

To do so, we reinforce the assumptions on the cone. Let us assume that the norm preserves the order relation $\prec$ over the cone $K$ (i.e., if $x,y\in K$ with $x\prec y$, then $\left\|x\right\|<\left\|y\right\|$).

	\begin{example}
		\begin{enumerate}[$(i)$]
			\item In $\mathbb{R}^2$, the usual Euclidean norm preserves the order relation $\prec$ given by the cone $K=\{(x_1,x_2)\in \mathbb{R}^2:x_1\geq 0, \ x_2\geq 0 \}$.
			\item In $\mathcal{C}([0,1])$, the standard sup-norm $\left\|\cdot\right\|_{\infty}$ does not preserve the order relation $\prec$ given by the cone of non-negative continuous functions $P=\{x\in \mathcal{C}([0,1]):x(t)\geq 0 \text{ for all } t\in [0,1] \}$. Indeed, the functions $y(t)=1$ and $x(t)=t$ for all $t\in [0,1]$ satisfy that $x,y\in P$ with $y-x\in P\setminus\{0\}$, but $\left\|y\right\|_{\infty}=\left\|x\right\|_{\infty}=1$.
			\item The sup-norm preserves the order relation $\prec$ given by the following sub-cone of $P$:
			\[K=\{x\in P:x \text{ is nonincreasing} \}. \]    
			Indeed, if $0\preceq x\prec y$, then there exists $h\in K\setminus\{0\}$ such that $y-x=h$. By definition of the cone $K$, $\left\|u\right\|_{\infty}=u(0)$ for every $u\in K$ and so $\left\|y\right\|_{\infty}=y(0)=x(0)+h(0)>x(0)=\left\|x\right\|_{\infty}$, since $h\neq 0$.
		\end{enumerate}
	\end{example}



We will need the following technical result.

\begin{lemma}
	Let $K$ be a cone in a normed linear space such that the norm preserves the order relation $\prec$ over $K$, $h\in K\setminus\{0\}$ and $r>0$. Then, for each fixed $x\in \overline{K}_r$, there exists $t_x\in [0,+\infty)$ satisfying that $\left\|x+t_x\,h\right\|=r$. Furthermore, the map 
	\begin{equation}\label{eq_map_t}
		\begin{array}{l}
			\boldsymbol{t}_{h}^{r}\equiv\boldsymbol{t}:\overline{K}_r\longrightarrow [0,+\infty) \\ \qquad \qquad x\longmapsto \boldsymbol{t}(x)=t_x
		\end{array}
	\end{equation}  
	is well--defined and continuous.
\end{lemma}

\noindent
{\bf Proof.} Given $x\in \overline{K}_r$, let us begin by showing the existence of $t_x$. Consider the continuous function $\phi:[0,\infty)\to\mathbb{R}$ defined as
$\phi(t)=\left\|x+t\,h\right\|.$
If $x\in K_r$, then $\phi(0)=\left\|x\right\|<r$ and $\lim_{t\to+\infty}\phi(t)=+\infty$, so Bolzano's theorem ensures the existence of $t_x>0$ such that $\phi(t_x)=r$. In case $x\in \partial K_r$, choose $t_x=0$.

Let us establish now the uniqueness of $t_x$. Suppose that there exist ${t}_x^1$ and ${t}_x^2$, with ${t}_x^1<{t}_x^2$, such that
\[\left\|x+{t}_x^1\,h\right\|=r=\left\|x+{t}_x^2\,h\right\|. \] 
Then $(x+{t}_x^2\,h)-(x+{t}_x^1\,h)=({t}_x^2-{t}_x^1)h\in K\setminus\{0\}$, that is, $x+{t}^1_x\,h\prec x+{t}_x^2\,h$. Since the norm preserves the order relation $\prec$ over $K$, we have $\left\|x+{t}_x^1\,h\right\|<\left\|x+{t}_x^2\,h\right\|$, a contradiction. Hence, the map $\boldsymbol{t}$ is well--defined.

Finally, let us prove that the map $\boldsymbol{t}$ is sequentially continuous. Let $\{x_n\}_{n\in\mathbb{N}}\subset \overline{K}_r$ be a convergent sequence to $\hat{x}\in \overline{K}_r$ and let us see that the sequence $\{\boldsymbol{t}(x_n)=t_{x_n} \}_{n\in\mathbb{N}}$ is convergent to $\boldsymbol{t}(\hat{x})=t_{\hat{x}}$. Since $\left\|x_n+t_{x_n}\,h\right\|=r$ for each $n\in\mathbb{N}$ and $\{x_n \}_{n\in\mathbb{N}}$ is convergent (so bounded), it follows that the sequence of real numbers $\{t_{x_n}\}_{n\in\mathbb{N}}$ is also bounded. Hence, it admits a convergent subsequence. Let $\{t_{x_{n_k}}\}_{k\in\mathbb{N}}$ be a convergent subsequence of $\{t_{x_n}\}_{n\in\mathbb{N}}$ with limit $\hat{t}$. By the continuity of the norm,
\[r=\lim_{k\to+\infty} \left\|x_{n_k}+t_{x_{n_k}}\,h\right\|=\left\|\hat{x}+\hat{t}\,h\right\|  \]
and thus the uniqueness ensures that $\hat{t}=t_{\hat{x}}$. Therefore, every convergent subsequence of $\{t_{x_n}\}_{n\in\mathbb{N}}$ has limit $t_{\hat{x}}$, which implies that the whole sequence $\{t_{x_n}\}_{n\in\mathbb{N}}$ converges to $t_{\hat{x}}$.
\qed

	\begin{remark}
		Note that the assumption that the norm preserves the order relation $\prec$ given by the cone can be replaced by the following: $h\in K\setminus\{0\}$ satisfies that
		\[x,y\in K \text{ such that } y-x=\lambda\,h \text{ for some } \lambda>0 \text{ imply } \left\|x\right\|<\left\|y\right\|. \]
	\end{remark}

	\begin{example}
		Consider the normed space $X=\{x\in\mathcal{C}([0,1]):x(0)=x(1)=0 \}$ endowed with the usual sup-norm and the cone $K=\{x\in X: x(t)\geq 0 \text{ for all } t\in [0,1] \}$. Then any function $h\in X$ such that $h(t)>0$ for all $t\in (0,1)$ satisfies the condition of the previous remark. 
	\end{example}	

As a consequence of the previous lemma, when the norm $\left\|\cdot\right\|_i$ preserves the order relation $\prec_i$ over $K_i$, the map $\hat{\rho}_i:(\overline{K}_i)_{R_i}\rightarrow(\overline{K}_i)_{r_i,R_i}$ defined as
\begin{equation}\label{eq_retract2}
	\hat{\rho}_i(x)=\left\{\begin{array}{ll} x+\boldsymbol{t}_{h_i}^{r_i}(x)h_i, & \text{ if } x\in (K_i)_{r_i}, \\ x, & \text{ if } x\in (\overline{K}_i)_{r_i,R_i}, \end{array} \right.
\end{equation}
with $h_i\in K_i\setminus\{0\}$ fixed, $i=1,2$, is a retraction.

\begin{lemma}
	Let $K_1$ and $K_2$ be two cones in the normed spaces $X$ and $Y$, respectively, such that the corresponding norms preserve the associated order relations $\prec_i$ over $K_i$, $i=1,2$. Let $T=(T_1,T_2):\overline{K}_{r,R}\rightarrow K$ be a compact map satisfying that for each $i\in\{1,2\}$ there exists $h_i\in K_i\setminus\{0\}$ such that the following condition holds in $\overline{K}_{r,R}$:
	\[x_i\neq T_i x+\mu\,h_i \ \text{ if } \left\|x_i\right\|=r_i \text{ and } \mu>0. \]
	
	If $T$ has no fixed point on the boundary of $K_{r,R}$, then its extension $N$ defined as $N:=T\circ \hat{\rho}$, with $\hat{\rho}(x):=(\hat{\rho}_1(x_1),\hat{\rho}_2(x_2))$ given by \eqref{eq_retract2}, has also no fixed point $x=(x_1,x_2)\in\overline{K}_{R}$ such that $\left\|x_i\right\|=r_i$ or $\left\|x_i\right\|=R_i$ for $i\in\{1,2\}$.		
\end{lemma}

\noindent
{\bf Proof.} Let us suppose that $N$ has a fixed point $(x_1,x_2)\in \overline{K}_{R}$ such that $\left\|x_1\right\|\in \{r_1,R_1\}$. Since $N=T$ in $\overline{K}_{r,R}$ and $T$ has no fixed points on the boundary of $\overline{K}_{r,R}$, it follows that $\left\|x_2\right\|<r_2$. Then
\[(x_1,x_2)=N(x_1,x_2)=\left(T_1(x_1,\hat{\rho}_2(x_2)),T_2(x_1,\hat{\rho}_2(x_2))\right).\]
In particular, $x_2=T_2(x_1,\hat{\rho}_2(x_2))$ and, by the definition of $\hat{\rho}_2$, we have $\hat{\rho}_2(x_2)=x_2+\boldsymbol{t}(x_2)h_2$. Therefore,
\[\hat{\rho}_2(x_2)-\boldsymbol{t}(x_2)h_2=T_2(x_1,\hat{\rho}_2(x_2)), \]
which is a contradiction since $\left\|\hat{\rho}_2(x_2)\right\|=r_2$ and $\boldsymbol{t}(x_2)>0$.
\qed

Now, repeating exactly the same fixed point index computations such as in the proof of Theorem~\ref{th_Korder} (which we omit here), it is possible to deduce the following fixed point theorem.

\begin{theorem}\label{th_Kras_hom}
	Let $r,R\in\R^2$, $r=(r_{1},r_{2})$, $R=(R_{1},R_{2})$, with $0<r_i<R_i$ $(i=1,2)$.	
	Assume that $K_1$ and $K_2$ are cones such that the corresponding norms preserve the order relations over $K_1$ and $K_2$, respectively, and $T=(T_1,T_2):\overline{K}_{r,R}\rightarrow K$ is a compact map such that for each $i\in\{1,2\}$ there exists $h_i\in K_i\setminus\{0\}$ satisfying the following condition in $\overline{K}_{r,R}$:
	\begin{enumerate}
		\item[$(a)$] $T_i x+\mu\,h_i\neq x_i$ if $\left\|x_i\right\|=r_i$ and $\mu>0$, and $T_i x\neq \lambda\, x_i$ if $\left\|x_i\right\|=R_i$ and $\lambda>1$.
	\end{enumerate}
	If $T$ has no fixed points on the boundary of $K_{r,R}$, then $i_{K}(T,K_{r,R})=1$. In particular, $T$ has at least one fixed point in $\overline{K}_{r,R}$.
\end{theorem}

\subsection{Multiplicity result}

The computation of the fixed point index is useful in order to obtain multiplicity results. In particular, it is possible to establish a three-solution criterion in the line of Amann \cite{amann}. For simplicity, we state the result with compression conditions of order type. Clearly, analogous conclusions can be derived under homotopy type conditions. 

\begin{theorem}
	\label{th_multi}
	Let $r^{(j)},R^{(j)}\in\R^2$, $r^{(j)}=(r^{(j)}_{1},r^{(j)}_{2})$, $R^{(j)}=(R^{(j)}_{1},R^{(j)}_{2})$, with $0<r^{(j)}_i<R^{(j)}_i$ $(i=1,2, \ j=1,2,3)$. Assume that the sets $\overline{K}_{r^{(j)},R^{(j)}}$ are such that 
	\[\overline{K}_{r^{(1)},R^{(1)}}\cup \overline{K}_{r^{(2)},R^{(2)}}\subset {K}_{r^{(3)},R^{(3)}} \quad \text{ and } \quad \overline{K}_{r^{(1)},R^{(1)}}\cap \overline{K}_{r^{(2)},R^{(2)}}=\emptyset. \]
	
	Moreover, assume that $T=(T_1,T_2):\overline{K}_{r^{(3)},R^{(3)}}\rightarrow K$ is a compact map and for each $i\in\{1,2\}$ and each $j\in\{1,2,3\}$ the following condition is satisfied in $\overline{K}_{r^{(j)},R^{(j)}}$:
	\begin{enumerate}
		\item[$(a^{\star})$] $T_i x\nprec x_i$  if $\left\|x_i\right\|=r^{(j)}_i$ and $T_i x\nsucc x_i$ if $\left\|x_i\right\|=R^{(j)}_i$.
	\end{enumerate}
	In addition, assume that $T$ is fixed point free on the boundaries of $\overline{K}_{r^{(1)},R^{(1)}}$ and $\overline{K}_{r^{(2)},R^{(2)}}$. Then $T$ has at least three distinct fixed points $\bar{x}^1$, $\bar{x}^2$ and $\bar{x}^3$ such that 
	\[\bar{x}^1\in K_{r^{(1)},R^{(1)}}, \ \ \bar{x}^2\in K_{r^{(2)},R^{(2)}} \ \ \text{ and } \ \ \bar{x}^3\in \overline{K}_{r^{(3)},R^{(3)}}\setminus\left(\overline{K}_{r^{(1)},R^{(1)}}\cup \overline{K}_{r^{(2)},R^{(2)}} \right). \]
\end{theorem}

\noindent
{\bf Proof.} As a consequence of Theorem \ref{th_Korder}, it follows that
\[i_{K}(T, {K}_{r^{(1)},R^{(1)}})=i_{K}(T, {K}_{r^{(2)},R^{(2)}})=1. \]
Then the existence property of the fixed point index ensures that $T$ has fixed points $\bar{x}^1\in K_{r^{(1)},R^{(1)}}$ and $\bar{x}^2\in K_{r^{(2)},R^{(2)}}$. Since $\overline{K}_{r^{(1)},R^{(1)}}\cap \overline{K}_{r^{(2)},R^{(2)}}=\emptyset$, they are two distinct fixed points.

Now, if $T$ has a fixed point, $\bar{x}^3$, on the boundary of $\overline{K}_{r^{(3)},R^{(3)}}$, then $\overline{K}_{r^{(1)},R^{(1)}}\cup \overline{K}_{r^{(2)},R^{(2)}}\subset {K}_{r^{(3)},R^{(3)}}$ implies that $\bar{x}^3$ is different from $\bar{x}^1$ and $\bar{x}^2$. Otherwise, if $T$ is fixed point free on the boundary of $\overline{K}_{r^{(3)},R^{(3)}}$, Theorem \ref{th_Korder} gives $i_{K}(T, {K}_{r^{(3)},R^{(3)}})=1$. By the additivity property of the index,
\[i_{K}(T, {K}_{r^{(3)},R^{(3)}}\setminus(\overline{K}_{r^{(1)},R^{(1)}}\cup \overline{K}_{r^{(2)},R^{(2)}}))=i_{K}(T, {K}_{r^{(3)},R^{(3)}})-i_{K}(T, {K}_{r^{(1)},R^{(1)}})-i_{K}(T, {K}_{r^{(2)},R^{(2)}})=-1, \]
which provides the third fixed point by means of the existence property of the index.
\qed

\begin{remark}
Observe that, if in Theorem \ref{th_multi} one requires $\overline{K}_{r^{(1)},R^{(1)}}\cup \overline{K}_{r^{(2)},R^{(2)}}\subset \overline{K}_{r^{(3)},R^{(3)}}$ (instead of $\overline{K}_{r^{(1)},R^{(1)}}\cup \overline{K}_{r^{(2)},R^{(2)}}\subset {K}_{r^{(3)},R^{(3)}}$) and that $T$ has not fixed points on $\partial_K\, K_{r^{(3)},R^{(3)}}$ then the conclusion holds for $\bar{x}^3\in K_{r^{(3)},R^{(3)}}\setminus\left(\overline{K}_{r^{(1)},R^{(1)}}\cup \overline{K}_{r^{(2)},R^{(2)}} \right)$.
\end{remark}

\section{Application to systems of boundary value problems}

Consider the following system of second-order differential equations subject to mixed boundary conditions
\begin{equation}\label{eq_BVP}
	\left\{\begin{array}{ll} x''(t)+f_1(t,x(t),y(t))=0 & \quad t\in [0,1], \\ y''(t)+f_2(t,x(t),y(t))=0 & \quad t\in [0,1], \\ x'(0)=x(1)=0=y'(0)=y(1), & \end{array} \right.
\end{equation}
where $f_1,f_2:[0,1]\times\mathbb{R}^{2}_{+}\rightarrow\mathbb{R}_{+}$ ($\mathbb{R}_+=[0,+\infty)$) are continuous functions satisfying that
\begin{enumerate}
	\item[$(\mathcal{H})$] for each $M_i>0$ there exists $C_i>0$ such that
	\[f_i(t,x_1,x_2)\leq C_i \ \text{ for all } x_i\in [0,M_i], \ x_j\in [0,+\infty) \ (j\neq i) \ \text{ and } \ t\in[0,1] \quad (i=1,2). \] 
\end{enumerate}

It is well--known that problem \eqref{eq_BVP} can be rewritten as an equivalent system of integral equations of the form
\begin{equation}\label{eq_sys}
	\left\{\begin{array}{l} x(t)=\displaystyle\int_{t}^{1}\int_{0}^{r} f_1(s,x(s),y(s))\,ds\, dr=:T_1(x,y)(t), \\ y(t)=\displaystyle\int_{t}^{1}\int_{0}^{r} f_2(s,x(s),y(s))\,ds\, dr=:T_2(x,y)(t). \end{array} \right. \quad (t\in [0,1])
\end{equation}
Let us consider the normed space of continuous functions $X=\mathcal{C}([0,1])$ endowed with the usual sup-norm $\left\|\cdot\right\|_{\infty}$ and the cones
\[K_1=K_2=\left\{u\in X:u(t)\geq 0 \text{ for all }t\in[0,1] \text{ and } u \text{ is nonincreasing on } [0,1] \right\}. \]
It is easy to check that the operator $T:=(T_1,T_2)$, defined as in \eqref{eq_sys}, maps the cone $K:=K_1\times K_2$ into itself and it is completely continuous. Simple calculations allow the operator $T$ to be expressed in terms of the corresponding Green's function in the following equivalent form: 
\[T_i(x,y)(t)=\displaystyle\int_{0}^{1}G(t,s)f_i(s,x(s),y(s))\,ds, \quad t\in [0,1] \quad (i=1,2), \]
where the kernel $G$ represents the Green's function of the second-order problem with Robin boundary conditions, which is given by
\[G(t,s)=\left\{\begin{array}{ll} 1-t & \text{ if } 0\leq s\leq t\leq 1, \\ 1-s & \text{ if } 0\leq t<s\leq 1, \end{array} \right. \]
see for instance \cite{WL}.

Now, we will provide eigenvalue criteria for the existence of positive (non-trivial in both components) solutions of the system \eqref{eq_BVP} following the ideas developed in \cite{WL}. To do so, we need to consider the compact linear Hammerstein integral operator
\begin{equation}\label{eq_Llinear}
	Lu(t):=\int_{0}^{1}G(t,s)u(s)\,ds, \quad t\in[0,1],
\end{equation}
and its corresponding eigenvalues, that is, the values $\lambda$ such that $\lambda u=Lu$. Let us denote by $\lambda_1$ the largest possible real eigenvalue of $L$ and by $\mu_1=1/\lambda_1$, the smallest positive characteristic value. 

\begin{remark}
	In the case of our linear operator $L$, finding an eigenvalue $\lambda$ of $L$ is equivalent to finding a non-trivial solution of 
	\[u''+\mu\,u=0, \quad u'(0)=u(1)=0, \]
	with $\mu=1/\lambda$. This linear differential problem only has nonzero solutions for $\mu>0$. In such case, solutions of the second--order differential equation are of the form:
	\[u(t)=A\cos\left(\sqrt{\mu }t\right)+B\sin\left(\sqrt{\mu}t\right), \quad A,B\in\mathbb{R}. \]
	Hence, the boundary condition $u'(0)=u(1)=0$ gives $B=0$ and $\mu=\left(\pi/2+k\,\pi\right)^2$, $k\in \mathbb{N}\cup\{0\}$. Therefore, $\mu_1=\pi^2/4$ and the associated eigenfunction is $h(t)=\cos\left(\pi t/2\right)$.
\end{remark}

Our existence result involves the value $\mu_1$ associated to the linear operator $L$ and the asymptotic behavior of the nonlinearities at $0$ and $+\infty$ expressed in terms of the following limits:
\[(f_i)_{0}:=\lim_{x_i\to 0^+}\dfrac{f_i(t,x_1,x_2)}{x_i} \quad \text{ and } \quad
(f_i)_{\infty}:=\lim_{x_i\to +\infty}\dfrac{f_i(t,x_1,x_2)}{x_i}. \]

\begin{theorem}\label{th_exis_BVP}
	Assume that for each $i\in\{1,2\}$ the following condition is satisfied:
	\[ \mu_1<(f_i)_{0}\leq+\infty \ \text{ and } \ 0\leq(f_i)_{\infty}<\mu_1, \text{ uniformly for }\ t,\  x_j \ (j\neq i). \]
	Then the system \eqref{eq_BVP} has at least one positive solution (non-trivial in both components).
\end{theorem}

\noindent
{\bf Proof.} For a fixed $i\in\{1,2\}$, let $\varepsilon>0$ be such that $(f_i)_{0}>\mu_1+\varepsilon$. Then there exists $\delta_i>0$ such that the following condition holds for each $x_i\in (0,\delta_i]$,
\begin{equation}\label{eq1_sup}
	f_i(t,x_1,x_2)\geq (\mu_1+\varepsilon)\,x_i \ \text{ uniformly for } t\in [0,1], \ x_j \ (j\neq i).
\end{equation} 
Fix $r_i\in (0,\delta_i)$ and let us prove that
\begin{equation}\label{eq_hom0}
	x_i\neq T_i x+\mu\,h \ \text{ if } \left\|x_i\right\|_{\infty}=r_i, \ x_j\in K_j \ (j\neq i) \text{ and } \mu>0,
\end{equation}
where $h$ is the eigenfunction of $L$ with $\left\|h\right\|_{\infty}=1$ associated to the eigenvalue $\lambda_1$.  

Let us suppose the contrary, that is, there exist $x\in K$ with $\left\|x_i\right\|_{\infty}=r_i$ and $\mu>0$ satisfying that $x_i= T_i x+\mu\,h$. Then
\[x_i(t)\geq \mu\,h(t) \ \text{ and } \ L x_i(t)\geq \mu\, L h(t)=\dfrac{\mu}{\mu_1}h(t), \quad t\in[0,1]. \]
Now, applying \eqref{eq1_sup}, we have
\begin{align*}
	x_i(t)&=\int_{0}^{1}G(t,s)f_i(s,x_1(s),x_2(s))\,ds+\mu\,h(t) \geq (\mu_1+\varepsilon) Lx_i(t)+\mu\,h(t) \\ &\geq (\mu_1+\varepsilon) \dfrac{\mu}{\mu_1}h(t)+\mu\,h(t) >2\mu\,h(t), \quad t\in[0,1].
\end{align*}
Repeating the reasoning, for each $n\in \mathbb{N}$ we obtain that $\left\|x_i\right\|_{\infty}>n\mu\left\|h\right\|_{\infty}=n\mu$, a contradiction.

On the other hand, let $\varepsilon>0$ ($\varepsilon<\mu_1$) be such that $(f_i)_{\infty}<\mu_1-\varepsilon$. Then there exists $M_i>r_i$ such that for each $x_i\geq M_i$ we have
\[f_i(t,x_1,x_2)\leq (\mu_1-\varepsilon)x_i \ \text{ for all } t, \ x_j \ (j\neq i). \]
By condition $(\mathcal{H})$, there exists $C_i>0$ such that 
\begin{equation}\label{eq1_sub}
f_i(t,x_1,x_2)\leq C_i+(\mu_1-\varepsilon)x_i \ \text{ for all } x_1,x_2\in[0,+\infty) \ \text{ and } \ t\in[0,1]. 
\end{equation} 
Since $1/\mu_1$ is the spectral radius of $L$, there exists $\left(\dfrac{1}{\mu_1-\varepsilon}I-L \right)^{-1}$, so we define 
\[R_0^i:=\left(\dfrac{1}{\mu_1-\varepsilon}I-L \right)^{-1}\left(\dfrac{C_i}{\mu_1-\varepsilon}\right).\]

Let us see that for each $R_i>R_0^i$ we have
\begin{equation}\label{eq_hom1}
	T_i x\neq \lambda\,x_i \ \text{ for all } \left\|x_i\right\|_{\infty}=R_i, \ x_j\in (\overline{K}_j)_{r_j,R_j} \ (j\neq i) \ \text{ and } \lambda>1.
\end{equation}
Reasoning by contradiction, we deduce from \eqref{eq1_sub} that
\[x_i(t)\leq \lambda\, x_i(t)=T_i x(t)=\int_{0}^{1}G(t,s)f_i(s,x_1(s),x_2(s))\,ds\leq (\mu_1-\varepsilon)L x_i(t)+C_i, \quad t\in[0,1], \]
which implies
\[\left(\dfrac{1}{\mu_1-\varepsilon}I-L \right)x_i(t)\leq \dfrac{C_i}{\mu_1-\varepsilon}, \quad t\in[0,1]. \]
By \cite[Lemma 4.1]{W}, we obtain 
\[x_i(t)\leq \left(\dfrac{1}{\mu_1-\varepsilon}I-L \right)^{-1}\left(\dfrac{C_i}{\mu_1-\varepsilon} \right), \quad t\in[0,1], \]
and thus $\left\|x_i\right\|_{\infty}\leq R_0^i<R_i$, a contradiction.

Therefore, $T$ satisfies conditions \eqref{eq_hom0} and \eqref{eq_hom1}. Consequently, Theorem~\ref{th_Kras_hom} ensures the existence of a fixed point of $T$, that is, a positive solution to the system \eqref{eq_BVP}. 
\qed

\begin{remark}
	Note that the proof of condition \eqref{eq_hom0} does not remain valid for the case $\mu=0$. Thus, even though $T$ satisfies condition $(a)$, it is unclear whether it satisfies $(a^{\dag})$, so the existence result cannot be derived in a similar way from Theorem~\ref{th_KP_strong}.
\end{remark}

Finally, we illustrate the existence criterion with a concrete example.

\begin{example}
	The following system of second-order differential equations with mixed boundary conditions
	\[\left\{\begin{array}{ll} x''(t)+\dfrac{t+1}{y(t)+1}+x(t) \arctan(x(t))=0 & \quad t\in [0,1], \\ y''(t)+\dfrac{t}{e^{-x(t)}+1}+\dfrac{1}{3}\sqrt{y(t)}=0 & \quad t\in [0,1], \\ x'(0)=x(1)=0=y'(0)=y(1), & \end{array} \right.	\]
	has at least one positive solution which is non-trivial in both components.
	
	Observe that the functions $f_1$ and $f_2$, defined as
	\[f_1(t,x,y)=\dfrac{t+1}{y+1}+x \arctan(x) \ \text{ and } \ f_2(t,x,y)=\dfrac{t}{e^{-x}+1}+\dfrac{1}{3}\sqrt{y}, \]
	satisfy that
	\[(f_1)_{0^+}=+\infty, \ (f_1)_{+\infty}=\dfrac{\pi}{2}<\dfrac{\pi^2}{4}=\mu_1 \ \text{ uniformly for } t,\ y \]
	and
	\[(f_2)_{0^+}=+\infty, \ (f_2)_{+\infty}=0 \ \text{ uniformly for } t,\ x. \]
	Therefore, Theorem \ref{th_exis_BVP} gives the conclusion.
\end{example}

\section*{Acknowledgements}

The authors have been partially supported by Ministerio de Ciencia y Tecnología (Spain), AEI and Feder, grant PID2020-113275GB-I00, and Xunta de Galicia, Spain, Project ED431C 2023/12.

\end{document}